\documentclass{IJFS}

\usepackage{amsfonts,amsthm,amsmath,amssymb,dsfont}
\usepackage{amsmath,    enumitem}
\usepackage{graphicx}
\usepackage{color}
\usepackage[all]{xy}
\usepackage{subfigure}
\usepackage{tikz}
\usetikzlibrary[arrows,shapes,graphs,decorations.pathmorphing,backgrounds,positioning,fit,petri]
\usepackage{ amscd,   eufrak}

 \usepackage{amssymb}
 \usepackage{amsmath}
 
 \usepackage{mathrsfs}
 \usepackage{amsfonts}
\usepackage{placeins}
 \usepackage{verbatim}

\begin{document}

\addtocounter{page}{0}

\titleijfs{A Credibility Approach on Fuzzy Slacks Based Measure
(SBM) DEA Model}

\author[1]{D. Mahla} 
\author[1]{S. Agarwal}
%\author[1]{M. Hamidi}
\affil[1]{Department of Mathematics, Birla Institute of Technology and Science, Pilani, India}
%\affil[2]{Department of Mathematics, Shahid Beheshti University, Tehran, Iran}

\emails{deepakmahlabits@gmail.com,  shivi@pilani.bits-pilani.ac.in}

\CorrespondAuthor{D. Mahla}

\oddPageHead{A Credibility Approach on Fuzzy Slacks Based Measure
(SBM) DEA Model}
\evenPageHead{D. Mahla, S. Agarwal}

\abstractijfs{
Data Envelopment Analysis (DEA) is a multi-criteria technique based on linear programming to deal with many real-life problems, mostly in nonprofit organizations. The slacks-based measure (SBM) model is one of the DEA model used to assess the relative efficiencies of decision-making units (DMUs). The SBM DEA model directly used input slacks and output slacks to determine the relative efficiency of DMUs. In order to deal with qualitative or uncertain data, a fuzzy SBM DEA model is used to assess the performance of DMUs in this study. The credibility measure approach, transform the fuzzy SBM DEA model into a crisp linear programming model at different credibility levels is used. The results came from the fuzzy DEA model are more rational to the real-world situation than the conventional DEA model. In the end, the data of Indian oil refineries is collected, and the efficiency behavior of the companies obtained by applying the proposed model for its numerical illustration.
}

\keywords_ijfs{Data envelopment analysis, Decision-making units, Fuzzy  slacks-based measure, Relative efficiency.
\vspace{1mm}
}

\Vol{*}	\No{*}	\Year{****}	\Pages{**}
\Received{**}	\Revised{**}		\Accepted{**}

{\let\newpage\relax\vspace*{.1mm}
\noindent\rule{\textwidth}{.3mm}\maketitle}

\thispagestyle{fancylogo}

\mabstract

\section{ Introduction}
Data Envelopment Analysis (DEA) is a robust non-parametric linear programming based technique to asses the performance of non-profit organizations. Many researchers used this technique in education, supply chain management, banking, transport, etc. The ranking of decision-making units (DMUs) using multiple conflicting criteria is one of the uses of DEA  \cite{J04}. There are different types of DEA models in literature, some of them are, Charnes, Cooper, Rhodes (CCR) \cite{CCR}, Banker, Charnes, Cooper (BCC) \cite{BCC}, slacks-based measure (SBM) \cite{Tone}, new slack model (NSM) \cite{Agarwal}, etc. The present study used the SBM DEA model, which was proposed by Tone \cite{Tone}, to evaluate the relative efficiencies of the decision-making units(DMUs). The SBM DEA model directly deals with the input excess and output shortfalls. 

DEA, as the name suggests, used frontier analysis to assess the relative efficiencies of DMUs where the minute change in data can significantly change the frontier. One of the challenges in real-world situations is that the available data might be present in an uncertain or qualitative form, or sometimes some data might be missing.  Therefore, the conventional DEA models are absurd to use with these types of data \cite{IJFS}. The fuzzy set theory, which was developed by Zadeh \cite{Zadeh} in 1965, integrates with DEA \cite{cooper2000data} to deal with this types of situations by creating more rational forms of DEA models. Many researchers applied DEA models to evaluate the efficiency of DMUs under fuzzy environments \cite{I18, A1}. Over the three-decade of research, the four primary approaches namely, the tolerance approach \cite{Wendell}, the $\alpha$-level-based approach \cite{Liu}, the fuzzy ranking approach \cite{Cheng}, the possibility approach \cite{possibility approach} used by researcher to solve the fuzzy DEA models. Wen et al. \cite{W13} extended the traditional DEA models to a fuzzy environment proposed a fuzzy DEA model based on credibility measure. Chen et al. \cite{C13} concluded that the use of the fuzzy SBM DEA model for estimating efficiency values not only represents the characteristic of the uncertainty of the efficiency values, it also presents the potential effect of risk volatility on efficiency values. Haiso et al. \cite{H11} also concluded that linguistic terms could not entirely fit with the conventional DEA models. Puri et al. \cite{P13} used fuzzy SBM DEA models to handle the imprecise data, and calculate the mix-efficiencies of State Bank of Patiala in the Punjab state of India. Wanke et al. \cite{W18} presented an analysis of the efficiency of Angolan banks using fuzzy DEA and stochastic DEA models based on the $\alpha-$ level approach and different tail dependence structure, respectively. Recently, Bakhtavar et al. \cite{B19} used a special risk prioritization algorithm by failure mode and effects analysis by SBM DEA model under fuzzy conditions. 
As per the Zadeh examination of the fuzzy theory, he concluded the possibility distribution of the fuzzy variable replicates the probability distribution of the random variable in probability theory. Fuzzy LP models can be considered the evolution of conventional LP models, in which fuzzy variables play the role of fuzzy coefficients, and fuzzy events construct the constraints controlling the model. The possibility theory determines the possibility of fuzzy events. Lertworasirikul et al. \cite{possibility approach, credibility approach} studied the fuzzy DEA models built by  Guo et al. \cite{G}, which took the possibility criterion and the necessity criterion as a measure and solved the ranking problem with two distinct approaches namely, the possibility approach and the credibility approach.
Recently, Agarwal \cite{Agarwal} applied possibility measure to solve the fuzzy SBM DEA model. The Possibility measure is used extensively, but it has no self-dual property, which is undoubtedly needed for practice. Liu and Liu \cite{Expected value} proposed credibility measure in 2002, which shows the self-dual character. The credibility theory which managed personal conviction degree numerically given by Liu \cite{ L2007} and refined it in his next researches \cite{Liu, L2010}. In the present study, the SBM DEA model is extended with a fuzzy environment for evaluating the efficiency of DMUs and solved by credibility measure. 

The rest of the paper organized as follows, section 2 recalls the basic SBM and fuzzy SBM DEA models. In section 3, the credibility measure is used to solve the fuzzy SBM DEA model. In section 4, the relative efficiency of Indian oil companies is calculated and compared with existing methods. In the end, the conclusion is given.

\section{Preliminaries}
This section discusses the SBM DEA model, the fuzzy numbers, and the fuzzy SBM DEA model in detail. Some definitions which will be used in this study are also discussed in this section.
\subsection{SBM DEA Model}
	  In 2001, Tone \cite{Tone} proposed the SBM DEA model, which measured the efficiency of DMUs by minimizing the ratio of the reduction rate input slacks to the expansion rate of output slacks. This model is units invariant in nature and monotone decreasing for input excess and output shortfall. Consider there are $m$-inputs, $n$-outputs, $r$-number of DMUs, ${x}_{iz}$= amount of i$^{th}$ input used by z$^{th}$ DMU, ${y}_{jz}$= amount of j$^{th}$ output used by z$^{th}$ DMU, $S^-_{iz}$= slack in the i$^{th}$ input of the z$^{th}$ DMU, $S^+_{jz}$= slack in the j$^{th}$ output of the z$^{th}$ DMU, and $\lambda_{oz}$ are intensity variables. Then, the SBM \cite{SSBM} DEA model for DMU$_z$ is given by, 
	\begin{align}
	\begin{split}	
	\text{min}\hspace{5mm} \rho_z =&\frac{t-\frac{1}{m}\sum_{i=1}^{m}S_{iz}^-{x}_{iz}}{t+\frac{1}{n}\sum_{j=1}^{n}S_{rz}^+{y}_{jz}}\\
	\text{subject to} \hspace{10mm} &\sum_{o=1}^{r}\lambda_{oz}{x}_{io}+S^-_{iz}=t{x}_{iz}  ~~~~~~\forall i=1,\cdots,m\\
	&	\sum_{o=1}^{r}\lambda_{oz}{y}_{jo}-S^+_{jz}=t{y}_{jz}~~~~~~\forall j=1,\cdots,n\\
	&\lambda_{oz} \geq 0,~S^-_{iz}\geq 0,~S^+_{jz}\geq 0,~t>0 ~~~~~~\forall o=1,~\cdots,r.
	\end{split}
	\end{align}
	The crisp SBM DEA model requires the accurate value of attributes. But, sometimes, in real-life situations, the available data is vague due to human error, qualitative in nature, or non-obtainable. The above problem is solved by using fuzzy set theory. The definition of the fuzzy number and fuzzy SBM DEA model is discussed in the next subsections. 
\subsection{Fuzzy Number} 
	The fuzzy set theory was introduced by Zadeh \cite{Zadeh} to deal with the vagueness  of  the data. A  fuzzy set on universal set $M$ is   defined by $\tilde{M}={(x,\mu_{\tilde{M}}(x))|x \in M ;~ \mu_{\tilde{M}}(x) ~\in~ [0
	,1]}$ in which $\mu_{\tilde{M}}(x)$ is called the membership function of the fuzzy set. The fuzzy numbers are special kind of fuzzy sets defined on real numbers $R$ which satisfy the following properties:
	\begin{enumerate}
	    \item Fuzzy numbers are normal  (i.e.$ \exists  x \in R ~: ~ \mu_{\tilde{M}}(x)= 1 $).
	    \item Fuzzy numbers are convex (i.e. $\mu_{\tilde{M}}(x) \geq \text{min} \{\mu_{\tilde{M}}(b), ~ \mu_{\tilde{M}}(a)\} ~\forall~ a \leq x \leq b)$.
	    
	\item  The  membership function of fuzzy numbers is an upper semi-continuous function.
	\end{enumerate}
    There are several types of fuzzy numbers; However, a triangular fuzzy number is used in this study. A {\bf{triangular fuzzy number $\tilde M$}} is a special kind of fuzzy number. It is fuzzy variable determined by triplet $(r,s,u)~\text{such that}~(r<s<u)$,  with membership function as,\\
\begin{equation} \label{eq1}
\begin{split}
 \mu(\tilde M) & = \frac{x-r}{s-r},~~ \text{if} \;  r\leq x \leq s \\
 & = \frac{u-x}{u-s}, ~~\text{if} \;  s\leq x \leq u \\
& = 0,~~~~~~~~~ \text{otherwise}.\\
\end{split}
\end{equation}

	\subsection{Fuzzy SBM DEA Model}
	SBM DEA model is changed into a fuzzy SBM DEA model by considering the inputs and outputs are fuzzy numbers.  The i$^{th}$ input of the z$^{th}$ DMU is indicated by $\tilde{x}_{iz}$, and the j$^{th}$ output of the z$^{th}$ DMU is indicated by $\tilde{y}_{jz}$, are fuzzy input and output for DMU${_z}$, respectively.
	\begin{align}
	\begin{split}	
	\text{min}\hspace{5mm} \rho_z =&\frac{t-\frac{1}{m}\sum_{i=1}^{m}S_{iz}^-/\tilde{x}_{iz}}{t+\frac{1}{n}\sum_{j=1}^{n}S_{rz}^+/\tilde{y}_{jz}}\\
	\text{subject to} \hspace{10mm} &\sum_{o=1}^{r}\lambda_{oz}\tilde{x}_{io}+S^-_{iz}=t\tilde{x}_{iz}  ~~~~~~\forall i=1,\cdots,m\\
	&	\sum_{o=1}^{r}\lambda_{oz}\tilde{y}_{jo}-S^+_{jz}=t\tilde{y}_{jz}~~~~~~\forall j=1,\cdots,n\\
	&\lambda_{oz} \geq 0,~S^-_{iz}\geq 0,~S^+_{jz}\geq 0,~t>0 ~~~~~~\forall o=1,\cdots,r.
	\end{split}
	\end{align}

The model (3) is a fuzzy SBM DEA model, which can not be solved directly. The credibility measure based approach is used to convert the fuzzy SBM DEA model into a crisp linear programming problem (lpp). The credibility measure is discussed in the next section.

\section{Credibility Measure}
In this study, the fuzzy SBM DEA model is approached by credibility measure, and the credibility measure is defined as,\\
\textbf{Credibility Measure:} Consider ${\xi}$ be a nonempty set with $P{\{\xi\}}$ be the power set of ${\xi}$. Liu and Liu \cite{Expected value} defined the credibility set function Cr$\{.\}$ as credibility measure if it holds the following conditions:
\begin{enumerate}
\item Cr$\{{\xi}\}=1$,
\item Cr$\{Y\}\leq$Cr$\{Z\}$ whenever $Y\subset Z \in \xi$,
\item Cr$\{Y\}+$Cr$\{Y\}^C=1$ for any event $Y \in \xi$,
\item Cr$\{\cup_{i}Y_{i}\}=$Sup$_i$Cr$\{Y_i\}$ for any events $Y_i$ with Sup$_i$Cr$\{Y_i\}<0.5$.
\end{enumerate}
The triplet $(\xi,P(\xi),Cr)$ are called the credibility space \cite{credibility theory}.
The theory of credibility of fuzzy events and chance-constrained programming (CCP) is used in this study to solve the fuzzy SBM DEA model. Wen et al. \cite{Wen} has given the following results which are used in solving process of our fuzzy model.\\
	 \textbf{Theorem 1.} Consider $\psi_1$ and $\psi_2$ are two fuzzy variables defined on credibility space $(\xi,P(\xi),Cr)$. If Cr$\{\psi_1=y \}$ and Cr$\{\psi_2=y \}$ are quasi concave, then 
	\begin{enumerate}
	\item  Cr$\{\psi_1+\psi_2\leq d \}\geq \alpha$ iff $(\psi_1)_{2(1-\alpha)}^U+(\psi_2)_{2(1-\alpha)}^U \leq d$,
	\item Cr$\{\psi_1+\psi_2\leq d \}\leq \alpha$ iff $(\psi_1)_{2(1-\alpha)}^U+(\psi_2)_{2(1-\alpha)}^U \geq d$. Here, $0.5 \leq \alpha \leq 1$.
	\end{enumerate}
 \textbf{Theorem 2.} Consider $(\psi)_{\alpha}^L$ and $(\psi)_{\alpha}^U$ are the lower and upper bounds of $\alpha$-cut of $\psi$, respectively. Then,
	\begin{enumerate}
	\item if $k \geq 0$, then $(k\psi)_{\alpha}^U=k(\psi)_{\alpha}^U$ and $(k\psi)_{\alpha}^L=k(\psi)_{\alpha}^L$,
	\item if $k \leq 0$, then $(k\psi)_{\alpha}^U=k(\psi)_{\alpha}^L$ and $(k\psi)_{\alpha}^L=k(\psi)_{\alpha}^U$.
	\end{enumerate}
The credibility distribution of triangular fuzzy number $\tilde M$ (2) is defined as,

\begin{equation} 
\begin{split}
  \text{Cr}{(\tilde M \leq b)} & = 0,~~~~~~~~~~~ \text{if} \;  r \geq b \\
 & = \frac{b-r}{2(s-r)},~~ \text{if} \;  r\leq b \leq s \\
 & = \frac{b-2s+u}{2(u-s)}, ~~\text{if} \;  s\leq b \leq u \\
& = 1,~~~~~~~~~ r \leq b.\\
\end{split}
\end{equation}

\begin{equation} 
\begin{split}
 \text{Cr}{(\tilde M \geq b)} & = 1,~~~~~~~~~~~ \text{if} \;  r \geq b \\
 & = \frac{2s-r-b}{2(s-r)},~~ \text{if} \;  r\leq b \leq s \\
 & = \frac{u-b}{2(u-s)}, ~~\text{if} \;  s\leq b \leq u \\
& = 0,~~~~~~~~~ r \leq b.\\
\end{split}
\end{equation}

According to credibility measure, converting fuzzy-chance constraints into their equivalent crisp ones \cite{1} in one particular confidence level $\alpha \geq 0.5 $ is as equation (6): 

\begin{equation} 
\begin{split}
 \text{Cr}{(\tilde M \leq b)}\geq \alpha & \iff (2-2 \alpha)s + (2 \alpha -1 ) u \\
 \text{Cr}{(\tilde M \geq b)}\geq \alpha & \iff (2-2 \alpha)s + (2 \alpha -1 ) r
\end{split}
\end{equation}
 \subsection{Credibility Approach on fuzzy SBM DEA model}
In the procedure to transform the model in the credibility programming SBM DEA model, each fuzzy coefficient is considered as a fuzzy variable, and each constraint is defined as a fuzzy event. The fuzzy SBM DEA model (3) becomes the following credibility SBM DEA model:  
	\begin{equation}
	\begin{split}
	\begin{cases}	
	\text{min}\hspace{5mm} f_z\\
	 \text{subject to}:\\
	 \text{Cr}\bigg\{\frac{t-\frac{1}{m}\sum_{i=1}^{m}S_{iz}^-/\tilde{x}_{iz}}{t+\frac{1}{n}\sum_{j=1}^{n}S_{rz}^+/\tilde{y}_{jz}}\leq f_z\bigg\} \geq  \alpha\\
	\text{Cr}\bigg\{\sum_{o=1}^{r}\lambda_{oz}\tilde{x}_{io}+S^-_{iz}-t\tilde{x}_{iz}=0\bigg\} \geq  \alpha ~~~~~~\forall i=1,\cdots,m\\
		\text{Cr}\bigg\{\sum_{o=1}^{r}\lambda_{oz}\tilde{y}_{jo}-S^+_{jz}-t\tilde{y}_{jz}=0\bigg\} \geq  \alpha ~~~~~~\forall j=1,\cdots,n\\
	~~~~~~\lambda_{oz} \geq 0,~~S^-_{iz}\geq 0,~~S^+_{jz}\geq 0,~~t>0 ~~~~~~\forall o=1,\cdots,r.
	\end{cases}
	\end{split}
	\end{equation}
	In order to solve the credibility SBM DEA model, the equality sign of constraints are converted into the inequality. Thus, model (7) becomes,\vspace{-0.36 cm}
	
	\begin{equation}
	\begin{split}
	\begin{cases}	
	\text{min}\hspace{5mm} f_z\\
	 \text{subject to}:\\\vspace{0.2 cm}
	~~~ \text{Cr}\bigg\{\frac{t-\frac{1}{m}\sum_{i=1}^{m}S_{iz}^-/\tilde{x}_{iz}}{t+\frac{1}{n}\sum_{j=1}^{n}S_{rz}^+/\tilde{y}_{jz}}\leq f_z\bigg\} \geq  \alpha\\ \vspace{0.2 cm}
	~~~ \text{Cr}\bigg\{\sum_{o=1}^{r}\lambda_{oz}\tilde{x}_{io}+S^-_{iz}- t\tilde{x}_{iz}\leq0 \bigg\} \geq  \alpha ~~~~~~\forall i=1,\cdots,m\\ \vspace{0.2 cm}
	~~~ \text{Cr}\bigg\{\sum_{o=1}^{r}\lambda_{oz}\tilde{x}_{io}+S^-_{iz} - t\tilde{x}_{iz}\geq 0 \bigg\} \geq  \alpha ~~~~~~\forall i=1,\cdots,m\\ \vspace{0.2 cm}
		~~~ \text{Cr}\bigg\{\sum_{o=1}^{r}\lambda_{oz}\tilde{y}_{jo}-S^+_{jz} - t\tilde{y}_{jz}\leq0\bigg\} \geq  \alpha ~~~~~~\forall j=1,\cdots,n\\ \vspace{0.2 cm}
		~~~ \text{Cr}\bigg\{\sum_{o=1}^{r}\lambda_{oz}\tilde{y}_{jo}-S^+_{jz} - t\tilde{y}_{jz}\geq0\bigg\} \geq  \alpha ~~~~~~\forall j=1,\cdots,n\\ \vspace{0.2 cm}
	~~~~~~\lambda_{oz} \geq 0,~~S^-_{iz}\geq 0,~~S^+_{jz}\geq 0,~~t>0 ~~~~~~\forall o=1,\cdots,r.
	\end{cases}
	\end{split}
	\end{equation}
If the membership functions of the fuzzy variable are normal and convex and $\alpha \geq 0.5$, then by using theorem 1 and 2, the credibility programming SBM model (8) is transformed into the model (9):
\begin{equation}
	 \begin{split}
	 \begin{cases}	
	\text{min}\hspace{5mm} f_z\\
	 \text{subject to}:\\ 
	 ~~~  \bigg\{\frac{t-\frac{1}{m}\sum_{i=1}^{m}S_{iz}^-/\tilde{x}_{iz}}{t+\frac{1}{n}\sum_{j=1}^{n}S_{rz}^+/\tilde{y}_{jz}}\bigg\}_{2(1-\alpha)}^L \leq f_z\\ \vspace{0.2 cm}
	~~~ \bigg\{\sum_{o=1}^{r}\lambda_{oz}\tilde{x}_{io}+S^-_{iz}- t\tilde{x}_{iz}\bigg\}_{2(1-\alpha)}^U\leq 0 ~~~~~~\forall i=1,\cdots,m\\ \vspace{0.2 cm}
	~~~ \bigg\{\sum_{o=1}^{r}\lambda_{oz}\tilde{x}_{io}+S^-_{iz} - t\tilde{x}_{iz}\bigg\} _{2(1-\alpha)}^L\geq 0 ~~~~~~\forall i=1,\cdots,m\\ \vspace{0.2 cm}
		~~~ \bigg\{\sum_{o=1}^{r}\lambda_{oz}\tilde{y}_{jo}-S^+_{jz} - t\tilde{y}_{jz}\bigg\} _{2(1-\alpha)}^U\leq 0 ~~~~~~\forall j=1,\cdots,n\\ \vspace{0.2 cm}
		~~~ \bigg\{\sum_{o=1}^{r}\lambda_{oz}\tilde{y}_{jo}-S^+_{jz} - t\tilde{y}_{jz}\bigg\} _{2(1-\alpha)}^L \geq 0~~~~~~\forall j=1,\cdots,n\\ \vspace{0.2 cm}
	~~~~~~\lambda_{oz} \geq 0,~~S^-_{iz}\geq 0,~~S^+_{jz}\geq 0,~~t>0 ~~~~~~\forall o=1,\cdots,r.
	\end{cases}
	\end{split}
	\end{equation} 
	
	The lower $\alpha$-cut of triangular fuzzy number $\frac{\tilde a}{\tilde b}$, given as, $\{\frac{\tilde a}{\tilde b}\}_{\alpha}^L={\tilde a}_\alpha ^L*1/({\tilde b}_\alpha ^U)$. Therefore, the model (9) is proportionate to the model (10) as follows:
	\begin{equation}
	\begin{split}
	\begin{cases}	
	~~~ \text{min}\hspace{5mm} \frac{\bigg\{{t-\frac{1}{m}\sum_{i=1}^{m}S_{iz}^-/\tilde{x}_{iz}}\bigg\}_{2(1-\alpha)}^L}{\bigg\{{t+\frac{1}{n}\sum_{j=1}^{n}S_{rz}^+/\tilde{y}_{jz}}\bigg\}_{2(1-\alpha)}^U }\\
	 \text{subject to}:\\ \vspace{0.2 cm}
	~~~ \bigg\{\sum_{o=1}^{r}\lambda_{oz}\tilde{x}_{io}+S^-_{iz}- t\tilde{x}_{iz}\bigg\}_{2(1-\alpha)}^U\leq 0 ~~~~~~\forall i=1,\cdots,m\\ \vspace{0.2 cm}
	~~~ \bigg\{\sum_{o=1}^{r}\lambda_{oz}\tilde{x}_{io}+S^-_{iz} - t\tilde{x}_{iz}\bigg\} _{2(1-\alpha)}^L\geq 0 ~~~~~~\forall i=1,\cdots,m\\ \vspace{0.2 cm}
		~~~ \bigg\{\sum_{o=1}^{r}\lambda_{oz}\tilde{y}_{jo}-S^+_{jz} - t\tilde{y}_{jz}\bigg\} _{2(1-\alpha)}^U\leq 0 ~~~~~~\forall j=1,\cdots,n\\ \vspace{0.2 cm}
		~~~ \bigg\{\sum_{o=1}^{r}\lambda_{oz}\tilde{y}_{jo}-S^+_{jz} - t\tilde{y}_{jz}\bigg\} _{2(1-\alpha)}^L \geq 0~~~~~~\forall j=1,\cdots,n\\ \vspace{0.2 cm}
	~~~~~~\lambda_{oz} \geq 0,~~S^-_{iz}\geq 0,~~S^+_{jz}\geq 0,~~t>0 ~~~~~~\forall o=1,\cdots,r.
	\end{cases}
	\end{split}
	\end{equation} 
	The fractional model (10) can be changed into a crisp linear programming model by normalization given as follows:
	\begin{equation}
	\begin{split}	
	\begin{cases}
	\text{min}\hspace{5mm} 
	  \bigg\{{t-\frac{1}{m}\sum_{i=1}^{m}S_{iz}^-/\tilde{x}_{iz}}\bigg\}_{2(1-\alpha)}^L \\
	  \text{subject to}:\\ \vspace{0.2 cm}
~~~ \bigg\{{t+\frac{1}{n}\sum_{j=1}^{n}S_{rz}^+/\tilde{y}_{jz}}\bigg\}_{2(1-\alpha)}^U =1\\	\vspace{0.2 cm} 
	~~~ \bigg\{\sum_{o=1}^{r}\lambda_{oz}\tilde{x}_{io}+S^-_{iz}- t\tilde{x}_{iz}\bigg\}_{2(1-\alpha)}^U\leq 0 ~~~~~~\forall i=1,\cdots,m\\ \vspace{0.2 cm}
	~~~ \bigg\{\sum_{o=1}^{r}\lambda_{oz}\tilde{x}_{io}+S^-_{iz} - t\tilde{x}_{iz}\bigg\} _{2(1-\alpha)}^L\geq 0 ~~~~~~\forall i=1,\cdots,m\\ \vspace{0.2 cm}
		~~~ \bigg\{\sum_{o=1}^{r}\lambda_{oz}\tilde{y}_{jo}-S^+_{jz} - t\tilde{y}_{jz}\bigg\} _{2(1-\alpha)}^U\leq 0 ~~~~~~\forall j=1,\cdots,n\\ \vspace{0.2 cm}
		~~~ \bigg\{\sum_{o=1}^{r}\lambda_{oz}\tilde{y}_{jo}-S^+_{jz} - t\tilde{y}_{jz}\bigg\} _{2(1-\alpha)}^L \geq 0~~~~~~\forall j=1,\cdots,n\\ \vspace{0.2 cm}
	~~~~~~\lambda_{oz} \geq 0,~~S^-_{iz}\geq 0,~~S^+_{jz}\geq 0,~~t>0 ~~~~~~\forall o=1,\cdots,r.
	\end{cases}
	\end{split}
	\end{equation} 
	The fuzzy SBM DEA model is transformed into a crisp lpp model (11) using equation (6) and can be solved using any software programs like MATLAB, PYTHON, LINGO, etc. 
	\section{Numerical Problem}
	A numerical example is presented to illustrate the proposed methodology. Here, the data of Indian oil refineries is collected for the financial year 2017-18, and the relative efficiencies are calculated by using the proposed methodology. The Indian oil sector is one of the largest and core areas which influences the other core sectors effectively. India is one of the biggest importers of oil, and its imports are increasing year by year. Compare to 4.56 million barrels of oil per day consumption by India in 2016 which is increased to 4.69 million barrels of oil per day in 2017 \cite{ibef}. \\
The total of 249.4 Million Metric Tons (MMT) oil installed in their refineries made it the second-largest refiner in Asia. In India, there are 10 (6 PSUs+ 2 Joint Venture+ 2 private limited) oil companies with 23 refineries and 249.4 Million Metric Tons (MMT) installed capacity for the financial year 2017-18. 
The data is collected from the company's annual reports, Ministry of Petroleum and Natural Gas \cite{M}, and Price waterhouse Coopers (PWC) \cite{PWC}  to evaluate the company's relative efficiency. 
In this study, the efficiency of oil companies is computed using four attributes, two inputs as capital expenditure and oil throughput and two outputs as revenue and Nelson Complexity Index (NCI). The details of  inputs and outputs are discussed below:
\subsection {Inputs and Outputs}
\begin{itemize}
    
\item \textbf{Capital expenditure (CE):} Capital expenditure is funds utilized by an organization to require all the physical resources needed for the company. The CE makes new projects or new investments. This financial outlay is made by organizations to keep up or increment the extent of their tasks. CE can be used for creating new buildings, repairing the equipment, or constructing the room. In our model, CE is used as crisp input.

\item \textbf{Oil Throughput Utilization:}
Crude throughput is the aggregate sum of unrefined that goes into a refinery before it comes out processed. Oil throughput utilization is given in the company's annual reports; therefore, it should be taken as a crisp input.
\item
 \textbf{Revenue:} The amount of money generated by a company from the selling of the product within a stipulated period is revenue. Every company announces its annual revenue in its yearly reports precisely; therefore, it should be taken as a crisp output.
\item \textbf{Nelson Complexity Index (NCI):} In India, it has been observed due to the proper utilization of capacity levels and high gross refining margins (GRMs), the NCI is improving from the last few years. The NCI compares the cost of upgrading units to the price of the crude distillation unit. The precise calculation of NCI requires repeated construction of the process units and consistent, standardized public reporting of the cost of such installations. In reality, the creation of crude distillation units are made of different sizes with the various technologies and the data regarding the value of such construction is reported could be of poor quality and not homogeneous. Due to preexisting limitations and heterogeneous data, one must be cautious while processing the data and should expect a high level of uncertainty. Therefore, to deal with this uncertainty, NCI is converted into a fuzzy number in this study. 
\end{itemize}

In our particular data, we divide two different types of companies based on scale size. The large companies are those which having employees more than 10000, and the small companies are which having employees number less than 10000. The companies BPCL, ONGC, HPCL, IOCL, and RIL, are large companies and CPCL, NRL, BORL, and NEL are small companies. The relative efficiency of small companies calculated separately as the comparison of the relative efficiency of small companies with large companies is not reasonable. Still, large companies can be compared to small companies. The data is given in table 1 as,
\begin{table}[h]

	\caption{Data of Indian Oil Companies}
	\centering
	%\resizebox{\textwidth}{%
	\begin{tabular}{p{2.5 cm}|p{2 cm}|p{2 cm}|p{1 cm}|p{1.2 cm}}
        \hline  
		{\bf{Company}}& {\bf {Capital Expenditure (Rs. in Cr)}} &\bf {Crude oil throughput (MMT)}&\bf {NCI}&\bf {Revenue (Rs. in Cr)}\\ \hline
		BPCL	&	8161	&	28.2	&	5.8	&	279312.70	\\	
	ONGC+MRPL	&	73664	&	16.27	&	9.5	&	928877	\\	
	HPCL	&	6722.45	&	28.7	&	12	&	243226.66	\\	
	CPCL	&	963	&	10.8	&	8.2	&	44227.2404	\\	
	IOCL	&	20345	&	69	&	9.5	&	506428	\\	
	NRL	&	387	&	2.8	&	10.5	&	15923.19	\\	
	BORL	&	3000	&	6.7	&	9.1	&	88454.4808	\\	
	RIL	&	40000	&	70.5	&	12.7	&	529120.00
	\\	
	NEL	&	961.66	&	20.7	&	11.8	&	86636.66	\\	\hline
	\end{tabular}%
%	}
	\end{table}	\\
The values of NCI have some uncertainty; therefore, in the proposed work, the NCI is converted into triangular fuzzy number by fuzzification process. The Saaty's scale is used to fuzzify NCI and fuzzified triangular fuzzy numbers are shown in table 2.
	\begin{table}[h]
	\centering
	\caption{Fuzzy data of Indian Oil Companies}
%		\resizebox{\textwidth}{%
	\begin{tabular}{p{2 cm}|p{2.5 cm}|p{1.5 cm}|p{2 cm}|p{1 cm}|p{1.2 cm}}
        \hline
		{\bf{Company Size}}&{\bf{Company }}& {\bf {Capital Expenditure}} &\bf {Crude oil throughput}&\bf {NCI}&\bf {Revenue}\\
			\hline \hline	
	&BPCL	&	1.89	&	4.20	&	(4,5,6)	&	3.41	\\	
	&ONGC+MRPL	&	9.00	&	2.85	&	(6,7,8)	&	9.00	\\	
	Large&HPCL	&	1.73	&	4.26	&	(7,8,9)	&	3.09	\\	
	&IOCL	&	3.21	&	8.83	&	(6,7,8)	&	5.36	\\	
	&RIL	&	5.34	&	9.00	&	(8,9,9)	&	5.56	\\	\hline \hline
	&&&&&\\
&CPCL	&	3.57	&	5.17	&	(6,7,8)	&	5.00	\\	
Small&NRL	&	2.03	&	2.08	&	(7,8,9)	&	2.44	\\	
&BORL	&	9.00	&	3.59	&	(6,7,8)	&	9.00	\\	
&NEL	&	3.56	&	9.00	&	(8,9,9)	&	8.84	\\	\hline

	\end{tabular}%
%	}
	\end{table}
Further, the efficiencies of all companies are calculated by the proposed approach on the fuzzy SBM DEA model, and calculations are done using MATLAB. The relative efficiencies of the companies using $\alpha-$ cut and possibility measure approach on SBM DEA are also computed. Table 3 shows the efficiencies of small and large companies using credibility measure as well as the ranking of the companies based on relative efficiency. The results of table 3 can be interpreted as, companies NRL, BORL, and NEL are efficient companies with relative efficiency 1. CPCL is inefficient among small companies at every credibility level from 0.5 to 1, with given four criteria. It is noticeable from the table that the value of relative efficiencies is increasing with the increase of credibility level of small companies. 
\begin{table}[h]
\caption{Relative efficiency of companies using credibility measure approach}
	\centering
%		\resizebox{\textwidth}{%
		\begin{tabular}{ c|c|c|c|c|c|c|c}
        \hline
	 	Credibility Level	&	0.5	&	0.6	&	0.7	&	0.8	&	0.9	&	1	&		\\ \hline \hline
		 		Large Company & & & & & & & Ranking\\ \hline
		BPCL	&	0.8571	&	0.8484	&	0.8387	&	0.8275	&	0.8148	&	0.8	&	3	\\
		ONGC+MRPL	&	1	&	1	&	1	&	1	&	1	&	1	&	1	\\
	HPCL	&	1	&	1	&	1	&	1	&	1	&	1	&	1	\\
		IOCL	&	0.7595	&	1	&	1	&	1	&	1	&	1	&	2	\\
		RIL	&	0.6655	&	0.6937	&	0.701	&	0.7014	&	0.7025	&	0.7029	&	4	\\ \hline \hline
 		Small Company & & & & & & & \\ \hline	
		CPCL	&	0.5748	&	0.6836	&	0.6854	&	0.6872	&	0.6391	&	0.6309	&	2	\\
	NRL	&	1	&	1	&	1	&	1	&	1	&	1	&	1	\\
		BORL	&	1	&	1	&	1	&	1	&	1	&	1	&	1	\\
		NEL	&	1	&	1	&	1	&	1	&	1	&	1	&	1	\\ \hline 
\end{tabular}
%}
\end{table}
It is also evident from Table 3 that ONGC and HPCL  are efficient among large companies at all credibility level. The BPCL and RIL are inefficient while BPCL is relatively more efficient than RIL at every credibility level with given four criteria.  The rank of IOCL is third among large companies followed by ONGC, HPCL (ranked 1), and BPCL at credibility level 0.5. The overall rank of IOCL is second, as it is efficient at every  credibility level except 0.5. Further, the efficiency of these companies are also computed using the "$\alpha-$cut approach" and "possibility approach" on fuzzy SBM DEA model. Tables 4 and 5 show the relative efficiency using a $\alpha-$cut approach for small and large companies.

   	  	\begin{table}[h]
\caption{Relative efficiency of small companies using $\alpha-$cut approach}
	\centering
%	\resizebox{\textwidth}{%
		\begin{tabular}{ c|c|c|c|c}
        \hline
        {\bf{Credibility level}}& {\bf {CPCL}} &\bf {NRL}&\bf {BORL}&\bf {NEL}\\
	 	\hline	\hline
	 	0.5	&	[0.5272,0.6415]	&	[1,1]	&	[1,1]	&	[1,1]	\\
	0.6	&	[0.5381,0.6295]	&	[1,1]	&	[1,1]	&	[1,1]	\\
	0.7	&	[0.5491,0.6176]	&	[1,1]	&	[1,1]	&	[1,1]	\\
	0.8	&	[0.5601,0.6057]	&	[1,1]	&	[1,1]	&	[1,1]	\\
	0.9	&	[0.5712,0.594]	&	[1,1]	&	[1,1]	&	[1,1]	\\
	1	&	[0.5823,0.5823]	&	[1,1]	&	[1,1]	&	[1,1]	\\ \hline \hline
	\bf {Rank}& 2&1&1&1\\ \hline 
\end{tabular}
%}
\end{table}
 	  	\begin{table}[h]
\caption{Relative efficiency of large companies using $\alpha-$cut approach}
	\centering
%	\resizebox{\textwidth}{%
		\begin{tabular}{ c|c|c|c|c|c}
        \hline
        {\bf{Credibility level}}& {\bf {BPCL}} &\bf {ONGC}&\bf {HPCL}&\bf {IOCL}&\bf {RIL}\\
	 	\hline	\hline	0.5	&	[0.5610,0.8750]	&	[1,1]	&	[0.6052,1]	&	[0.5512,1]	&	[0.4880,0.8643]	\\
	0.6	&	[0.6077,0.8718]	&	[1,1]	&	[0.6394,1]	&	[0.5869,1]	&	[0.5176,0.8153]	\\
	0.7	&	[0.6586,0.8684]	&	[1,1]	&	[0.6925,1]	&	[0.6258,1]	&	[0.5498,0.7743]	\\
	0.8	&	[0.7140,0.8649]	&	[1,1]	&	[0.7508,1]	&	[0.6681,0.9633]	&	[0.5848,0.7375]	\\
	0.9	&	[0.7747,0.8611]	&	[1,1]	&	[0.8651,1]	&	[0.7145,0.8332]	&	[0.6233,0.6970]	\\
	1	&	[0.8571,0.8571]	&	[1,1]	&	[1,1]	&	[[0.7596,0.7596]	&	[0.6655,0.6655]	\\ \hline
 \hline
	\bf {Rank}& 3&1&1&2&4\\ \hline 
\end{tabular}
%}
\end{table}
Table 6 shows the relative efficiency using the possibility approach for small and large companies.
\begin{table}[h]
\caption{Relative efficiency of companies using possibility measure approach}
%	\resizebox{\textwidth}{%
	\centering
		\begin{tabular}{ c|c|c|c|c|c|c|c}
        \hline 
		
	Credibility Level	&	0.5	&	0.6	&	0.7	&	0.8	&	0.9	&	1	&		\\ \hline
	\hline
	Large Company & & & & & & &Ranking \\ \hline
	BPCL	&	0.8333	&	0.8387	&	0.8437	&	0.8484	&	0.8529	&	0.8571	&	3	\\
	ONGC+MRPL	&	1	&	1	&	1	&	1	&	1	&	1	&	1	\\
	HPCL	&	1	&	1	&	1	&	1	&	1	&	1	&	1	\\
	IOCL	&	1	&	1	&	1	&	1	&	1	&	0.7595	&	2	\\
	RIL	&	0.7011	&	0.701	&	0.6983	&	0.6937	&	0.6892	&	0.6655	&	4	\\
\hline \hline
Small Company & & & & & & & \\ \hline
	CPCL	&	0.6358	&	0.635	&	0.6342	&	0.6335	&	0.6327	&	0.5415	&	2	\\
	NRL	&	1	&	1	&	1	&	1	&	1	&	1	&	1	\\
	BORL	&	1	&	1	&	1	&	1	&	1	&	1	&	1	\\
	NEL	&	1	&	1	&	1	&	1	&	1	&	1	&	1	\\
\hline
\end{tabular}
%}
\end{table}
The efficiency values are comparable for two approaches with our proposed approach when the value of $\alpha$ is moving from 0.5 to 1.0.
\section{Conclusions}
This present study has established a credibility approach to solve the fuzzy SBM DEA model. This approach converted the fuzzy model into a credibility SBM DEA model. The credibility measure, which is the mean of necessity and possibility measure and self-dual in nature used to find the desired level of confidence. The credibility approach deals with the fuzzy environment in the fuzzy SBM DEA model significantly. The method is more reasonable and conceivable as compares to the tolerance approach, $\alpha-$cut approach, fuzzy ranking approach, and possibility approach. In the numerical example, the relative efficiency and ranking of Indian oil and refineries are computed, and results are compared with existing methods to solve the fuzzy SBM DEA model.
\section*{Future Work} The present work proposed a new approach to solving the fuzzy SBM DEA model, and the numerical method demonstrated the compatibility of the proposed method. From the mathematical problem, we can conclude that there are two limitations. First, we can not rank the efficient DMUs directly from this method, and the other is, the targets can not be computed because the data is rationalized. In the future, we will extend this model to rank the efficient DMUs, and to find out the input and output targets. 
  \section*{Acknowledgement} The authors thank DST for providing partial support under FIST grant SR/FST/MSI-090/2013(C) and CSIR for providing financial assistance under JRF programme (Award No.: 09/719/(0079)/2017-EMR-{I}. The authors wish to express their appreciation for several excellent suggestions for improvements in this paper made by the referees.\\

\end{document}